\documentclass[10pt,journal,final,a4paper,oneside,twocolumn]{IEEEtran}

\pdfoutput=1

\usepackage{graphicx}
\usepackage{amsmath}
\usepackage{amsfonts}
\usepackage{amssymb}
\usepackage{subfigure}
\usepackage{cite}
\usepackage{color}
\usepackage{upgreek}

\newcommand{\grad}{\nabla }
\renewcommand{\div}{\nabla \cdot}

\newcommand{\vect}[1]{\boldsymbol{#1}}
\newcommand{\bc}{\begin{color}{red}}
\newcommand{\ec}{\end{color}}

\newcommand{\figref}[1]{Figure~\ref{#1}}
\newcommand{\secref}[1]{Section~\ref{#1}}
\newcommand{\subsecref}[1]{Section~\ref{#1}}
\newcommand{\tabref}[1]{Table~\ref{#1}}
\newcommand{\EE}[1]{\,10^{#1}}
\DeclareMathOperator{\jj}{j}
\DeclareMathOperator{\dd}{d}
\newcommand{\tocheck}[1]{{#1}}

\begin{document}

\title{Hierarchical bases preconditioner to enhance convergence of the CFIE with multiscale meshes}

\author{M. Righero, I. M. Bulai, M. A. Francavilla, F. Vipiana, \textit{Senior Member, IEEE}, M. Bercigli, A. Mori, M. Bandinelli, and G. Vecchi, \textit{Fellow,
IEEE} 
%
%
\thanks{M. Righero and M. A. Francavilla are with the Antenna and EMC Lab (LACE),
Istituto Superiore Mario Boella, Torino 10138, Italy (e-mail: righero@ismb.it, francavilla@ismb.it).}
\thanks{I. M. Bulai is with the Department of Mathematics, Universit\`a di Torino, Torino 10123, Italy (e-mail: iuliamartina.bulai@unito.it).}
\thanks{M. Bercigli., A. Mori, and M. Bandinelli are with Ingegneria Dei Sistemi (IDS), Pisa, Italy (email m.bercigli@idscorporation.com, a.mori@idscorporation.com, m.bandinelli@idscorporation.com)}
\thanks{
F. Vipiana and G. Vecchi are with the Antenna and EMC Lab (LACE),
Politecnico di Torino, Torino 10129, Italy (e-mail: francesca.vipiana@polito.it, giuseppe.vecchi@polito.it).}
}

\maketitle

\begin{abstract}
A hierarchical quasi-Helmholtz decomposition, originally developed to address the low-frequency and dense-discretization breakdowns for the EFIE, is applied together with an algebraic preconditioner to improve the convergence of the CFIE in multiscale problems. The effectiveness of the proposed method is studied first on some simple examples; next, test on real-life cases up to several hundreds wavelengths show its good performance. 
\end{abstract}


\section{Introduction} \label{sec:introduction}

The Electric Field Integral Equation (EFIE) is successfully used to numerically simulate the electromagnetic behaviour of structures in the frequency domain. Its Galerkin discretization using the classical Rao-Wilton-Glisson (RWG) functions defined in \cite{1982_Rao_ETAL_Electromagnetic}, however, produces an ill-conditioned matrix when the frequency is low or the discretization density is high. A hierarchical quasi-Helmholtz decomposition can help overcoming this issue, recasting the linear system in a form that can be easily preconditioned with a diagonal rescaling and an algebraic preconditioner on the fewer functions at the coarsest level. 
The EFIE suffers from other convergence problems, mainly related to resonances in closed bodies. A possible cure is to use the Combined Field Integral Equation (CFIE) instead, which mixes the EFIE with the Magnetic Field Integral Equation (MFIE), and exploits the fact that resonances occur at different frequencies for the two formulations (e.g. \cite{2008_Gibson_TheMethod}). 
One is then tempted to use both the approaches and to apply the hierarchical quasi-Helmholtz decomposition to the CFIE, hoping to gain both the improvements. Focusing on the mesh density issue, we show that this can be done indeed, with simple examples and realistic cases to illustrate its effectiveness. Preliminary investigations along these lines were presented in the conference paper \cite{2014_FrancavillaETAL_Applications}.

Finally, we observe that the work in \cite{2015_AdrianETAL_Hierarchical, 2015_AdrianETAL_OnTheHierarchical} shares our endeavour, yet with a completely different approach. 


\section{Background}\label{sec:background}

The scattered electric and magnetic fields due to a time-harmonic surface current distribution on $\partial \Omega$, the boundary of the region $\Omega$, in a medium with permittivity $\epsilon$ and permeability $\mu$ are given by
\begin{equation}
\label{eq:EFIO}
\begin{split}
\vect{e}^s(\vect{r}) & = -\jj\omega\mu\int_{\partial \Omega}g(\vect{r}, \vect{r}')\vect{j}(r')\dd \vect{r}' \\ 
& + \frac {1}{\jj \omega \epsilon}\grad \int_{\partial \Omega} g(\vect{r}, \vect{r}') \div \vect{j}(r')\dd \vect{r}'
\end{split}
\end{equation}
and
\begin{equation}
\label{eq:MFIO}
\vect{h}^s(\vect{r}) = \int_{\partial \Omega}\grad g(\vect{r}, \vect{r}') \times \vect{j}(r')\dd \vect{r}' 
\end{equation}
respectively, where the neglected time factor is $\exp(\jj\omega t)$ and $g(\vect{r}, \vect{r}')=\exp(-\jj k |\vect{r}-\vect{r}'|)/4\uppi|\vect{r}-\vect{r}'|$ is the scalar Green's function relevant to the surrounding medium, whose wavenumber is $k$.
When considering a perfect electric conductor with impinging time-harmonic electric and magnetic fields $\vect{e}^i$ and $\vect{h}^i$, the unknown induced current $\vect{j}$ can be determined writing an equation on the boundary $\partial \Omega$, the EFIE
\begin{equation}
\label{eq:EFIE}
-\vect{n}(\vect{r}) \times \vect{n}(\vect{r}) \times ( \vect{e}^s(\vect{r}) + \vect{e}^i(\vect{r})) = \vect{0}, \vect{r}\in\partial \Omega.
\end{equation}
If the region is closed, the unknown induced current $\vect{j}$ can de determined writing an equation as approaching the boundary $\partial \Omega$ from the exterior, the MFIE
\begin{equation}
\label{eq:MFIE}
-\vect{n}(\vect{r}) \times ( \vect{h}^s(\vect{r}) + \vect{h}^i(\vect{r})) = \vect{j}(\vect{r}), \vect{r}\rightarrow\partial \Omega.
\end{equation}
Expanding the unknown current in a finite basis, as the RWG defined in \cite{1982_Rao_ETAL_Electromagnetic}, and testing \eqref{eq:EFIE} and \eqref{eq:MFIE} against the same functions, in a Galerkin or Method of Moments (MoM) fashion (e.g. \cite{2008_Gibson_TheMethod}), we obtain two systems of linear equations
\begin{equation}
\label{eq:Galerkin}
[Z][j]=[e^i], \quad [B][j]=[h^i]
\end{equation}
%
Solving either one  of \eqref{eq:Galerkin} with a direct or iterative method gives the searched coefficients. To have better convergence properties, especially when considering closed objects, a linear combination of the two systems in \eqref{eq:Galerkin} is considered
\begin{equation}
\label{eq:CFIE_Galerkin}
\left( 1/\eta[Z] + \alpha[B] \right) [j]=\left( 1/\eta[e^i] + \alpha[h^i] \right)
\end{equation}
called the CFIE, where $\eta=\sqrt{\epsilon\mu}$ and generally $\alpha=1/2$ or $\alpha=1$ in the literature (e.g. \cite{2008_Gibson_TheMethod}).
\section{CFIE with hierarchical bases}\label{sec:cfie_hierarchical}
The EFIE \eqref{eq:EFIE} produces ill-conditioned matrix $[Z]$ when applied at low frequencies or on fine discretizations. This is due to the presence of the divergence operator and the opposite scaling of the two parts in \eqref{eq:EFIO}, \cite{2010_AnriulliETAL_Solving}. A possible cure to this problem is based on Calderon's identities \cite{2008_Andriulli_ETAL_AMultiplicative}.  
Another options is to use a hierarchical change of basis, with functions with increasing support. 
As described in \cite{2008_Andriulli_ETAL_Hierarchical, 2009_Vipiana_Vecchi_ANovel}, this can be accomplished with a pre- and post-multiplication of the matrix $[Z]$, obtained in standard RWG functions, by a matrix $[Q]$ and its transpose, with a complexity scaling as $N \log N$, where $N$ is the number of unknowns. We refer to the hierarchical basis as Multi Resolution (MR) basis.
%
%
This  preconditioning effect is extremely efficient for very dense meshes, and it becomes less effective when the size of the support of the coarser functions is comparable with $\lambda/8$, $\lambda$ being the wavelength of the surrounding medium. As described in \cite{2010_Vipiana_ETAL_EFIEModeling}, the method is applied to large multiscale problems by: a) stopping the hierarchical construction at the desired level, and b) applying an algebraic preconditioner---as the LU or Incomplete LU (ILU) factorization---only on the part corresponding to the coarsest functions; that is precisely the regime where algebraic preconditioners are effective in \emph{smooth} structures.
The use of hierarchical loops as well as of hierarchical non-solenoidal functions is necessary to stop the hierarchical construction at a given level, and generate \emph{generalized} RWG functions \cite{2008_Andriulli_ETAL_Hierarchical} at the chosen coarse level; this, and the ability to efficiently find and use (possible) global loops constitute the dominant reasons for the use of hierarchical loops.

%

From an algorithmic point of view, the MR scheme can be directly applied to the CFIE. It amounts to a pre- and post-multiplication of the Galerkin matrix in \eqref{eq:CFIE_Galerkin} and a subsequent algebraic preconditioner, oblivious of the considered operator. 

However, the effectiveness of this approach is far from obvious.
%
%
%
In fact, the identity operator appearing in \eqref{eq:MFIE} questions the utility of using the MR basis with the CFIE. Indeed, while the MR change of basis effectively preconditions the EFIE, it could have a negative effect on the spectral properties of the operator of the MFIE, as explained in \cite{2015_AdrianETAL_Hierarchical, 2015_AdrianETAL_OnTheHierarchical}. From this perspective, the first question is if there is a trade-off between the regularizing effect on the EFIE part and the possible deterioration of the MFIE spectrum.

As our results will show, in all investigated cases, including electrically large and complex structures, the use of the proposed hybrid hierarchical approach is \emph{always} advantageous. We remark that the spectral properties of the CFIE are not simply the sum of the spectral properties of the constituent EFIE and MFIE;  we also remark that the CFIE, in practice, only exists at the level of (finite) matrices and coefficients. A full theoretical explanation of this findings is beyond the scopes of this communication; the spectral properties of a similar form of the CFIE are studied, inter alia, in \cite{2015_SimonETAL_OnTheHierarchical}; we are not aware of any other similar work on this topic.
%
%
%
%

\section{Numerical examples} \label{sec:examples}


In the first place, we verify that the MR basis does not corrupt the resonance free properties of the CFIE, considering an example with a sphere of radius $0.5$ m illuminated by a linearly polarized plane wave at different frequencies. In particular, we hit the resonant frequency of $474.56$ MHz. \tabref{tab:resonances} displays the number of GMRES iterations needed to reach a relative residual lower than $\EE{-4}$. The CFIE formulation does not suffer around the resonant frequency and the MR change of basis preserves this property (as it was reasonable to be expected, as it does not alter the constituent equations) and actually improves the convergence rate.
\begin{table}
\def\U{\rule{0pt}{3ex}}
\def\D{\rule[-0.5ex]{0pt}{0pt}}
\def\V{\U\D}
\begin{center}
\begin{tabular}{ccccc}
\hline
\V  & \multicolumn{2}{c}{EFIE} & \multicolumn{2}{c}{CFIE}\\
 \cline{2-3} \cline{4-5}
\V frequency [MHz] & RWG & MR & RWG &  MR \\
\hline
\U $300 $		& 250 & 275 & 76 & 57 \\
$400 $ 		& 278 & 303 & 72 & 57 \\
$474.56 $ 	& 466 & 470 & 70 & 60 \\
$500 $		& 250 & 408 & 70 & 56 \\
\D $600 $		& 250 & 540 & 69 & 56 \\
\hline
\end{tabular}
\end{center}
\caption{Number of GMRES iterations to reach a residual lower than $\EE{-6}$ for the resonant sphere test case, \secref{sec:examples}.}
\label{tab:resonances}
\end{table}

\subsection{Role of hierarchical loops} \label{subsec:hierarchical}

We focus on the preconditioning effect of the MR basis on complex structures with heterogeneous mesh and fine details; as reported in \secref{sec:introduction},  an algebraic preconditioner on the coarsest level complements the MR preconditioning. Inter alia, this implies use of hierarchical loops. 
We start with the analysis of the dense-discretization breakdown, as a function of the edge refinement. We consider a cube of edge $1$ m, at the frequency of $1$ MHz, and meshes with decreasing edge length, denoted with $h$ in \tabref{tab:hierarchical}. We investigate the condition number and the number of GMRES iterations to reach a residual lower than $\EE{-4}$, for an impinging linearly polarized plane wave. \tabref{tab:hierarchical} reports the results for the EFIE and the CFIE, both in the RWG, in the MR basis built with standard (point) loops, and in the MR basis built with hierarchical loops \cite{2009_Vipiana_Vecchi_ANovel}. Consistent with results in \cite{2009_Vipiana_Vecchi_ANovel}, there is a slight improvement with hierarchical loops for the EFIE (with removal of a single very small eigenvalue). 
We note instead that the effect of hierarchical loops is significantly larger for the CFIE than for the EFIE. 

%
\begin{table*}
\def\U{\rule{0pt}{3ex}}
\def\D{\rule[-0.5ex]{0pt}{0pt}}
\def\V{\U\D}
\begin{center}
\begin{tabular}{c p{.1mm} ccc p{.1mm}ccc }
\hline
\V  && \multicolumn{3}{c}{EFIE} && \multicolumn{3}{c}{CFIE} \\
\cline{3-5} \cline{7-9}
\V $1/h$ && RWG & MR hierarchical loops & MR point loops &&  RWG & MR hierarchical loops & MR point loops \\
\hline
\U 1&& $3.40 \EE{4}$ (14) & $5.51 \EE{0}$  (9) & $1.12 \EE{1}$  (9) && $1.25 \EE{2}$  (16) & $5.25 \EE{0}$  (8) & $1.39 \EE{1}$  (10) \\
2  	&& $1.42 \EE{5}$ (34) & $1.88 \EE{1}$ (15) & $5.86 \EE{1}$ (14) && $3.08 \EE{2}$ (33) & $1.20 \EE{1}$  (15) & $7.55 \EE{1}$ (18) \\
4  	&& $5.23 \EE{5}$ (51) & $3.92 \EE{1}$ (21) & $3.37 \EE{2}$ (21) && $7.43 \EE{2}$ (54) & $2.24 \EE{1}$  (20) & $4.31 \EE{2}$ (33) \\
8  	&& $2.07 \EE{6}$ (74) & $8.09 \EE{1}$ (27) & $1.74 \EE{3}$ (30) && $1.66 \EE{3}$ (84) & $4.16 \EE{1}$ (25) & $2.47 \EE{3}$ (70) \\
16 	&& $8.29 \EE{6}$ (107) & $1.93 \EE{2}$ (37) & $6.09 \EE{3}$ (43) && $3.56 \EE{3}$ (129) & $7.34 \EE{1}$ (35) & $1.30 \EE{4}$ (126) \\
20 	&& $1.39 \EE{7}$ (74) & $2.25 \EE{2}$ (41) & $1.94 \EE{4}$ (48) && $4.10 \EE{3}$ (119) & $6.54 \EE{1}$ (33) & $1.80 \EE{4}$ (157) \\
\hline
\end{tabular}
\end{center}
\caption{Condition number (and GMRES iterations to reach a residual lower than $\EE{-4}$) for the cube example, \subsecref{subsec:hierarchical}.}
\label{tab:hierarchical}
\end{table*}

\subsection{Role of the MR change of basis to handle details} \label{subsec:details}

We apply the MR change of basis and the preconditioning strategy devised in \cite{2010_Vipiana_ETAL_EFIEModeling} to a closed object with a heterogeneous mesh with fine details. We consider a sphere of radius $0.5$ m at the frequency of $300$ MHz. Starting from a regular mesh with edge length of approximately $\lambda/10$, we refine the mesh on part of the structure to have edges as short as $\lambda/200$. We investigate the convergence as the level of detail increases. In \tabref{tab:details}, we compare the number of GMRES iterations to reach a residual lower than $\EE{-4}$, when the excitation is a linearly polarized plane wave. We use the standard RWG basis (RWG) and the MR basis with an algebraic preconditioner on the coarsest level (MR), both for the EFIE and the CFIE. The maximum (minimum) area of the mesh cells, $A^{\mathrm{max}}$ ($A^{\mathrm{min}}$), and the maximum (minimum) length of the mesh cells edges, $h^{\mathrm{max}}$ ($h^{\mathrm{min}}$), keep trace of the increasing heterogeneity of the mesh. Even in this simple case, we can appreciate the better convergence of the CFIE coupled with the MR and the algebraic preconditioner.
\begin{table}
\def\U{\rule{0pt}{3ex}}
\def\D{\rule[-0.5ex]{0pt}{0pt}}
\def\V{\U\D}
\begin{center}
\begin{tabular}{cccccc}
\hline
\V &&  \multicolumn{2}{c}{EFIE} & \multicolumn{2}{c}{CFIE}\\
\cline{3-4} \cline{5-6}
\V $A^{\mathrm{max}}/A^{\mathrm{min}}$ & $h^{\mathrm{max}}/h^{\mathrm{min}}$ & RWG & MR & RWG &  MR \\
\hline
\U $9.44 \EE{0}$		& $4.23 \EE{0}$ &  $93$ & $21$ & $25$ & $21$ \\
$2.36 \EE{2}$ 		 	& $1.98 \EE{1}$ & $186$ & $39$ & $45$ & $28$ \\
\D $1.00 \EE{3}$ 		& $3.75 \EE{1}$ & $212$ & $53$ & $59$ & $35$ \\
\hline
\end{tabular}
\end{center}
\caption{Number of GMRES iterations to reach a residual lower than $\EE{-4}$ for the example of the heterogeneous mesh, \subsecref{subsec:details}.}
\label{tab:details}
\end{table}

\section{Real-life cases} \label{sec:realistic}

Realistic cases are simulated with a MoM code implementing the MultiLevel Fast Multipole Algorithm (e.g. \cite{2008_Gibson_TheMethod}) to evaluate the Galerkin testing (ADF-EMS \cite{ADF_web}). We compare the results obtained with an iterative solver in the RWG basis with a simple diagonal preconditioner (RWG), in the RWG basis with a full near-field LU preconditioner (RWG LU), and in the MR basis with a LU preconditioner on the near-field interactions relevant to the coarsest functions \cite{2010_Vipiana_ETAL_EFIEModeling} (MR), both enforcing the EFIE and the CFIE on the whole structure.

%
\begin{table*}
\def\U{\rule{0pt}{3ex}}
\def\D{\rule[-0.5ex]{0pt}{0pt}}
\def\V{\U\D}
\begin{center}
\begin{tabular}{l p{.1mm} ccc p{.1mm}ccc }
\hline
\V  && \multicolumn{3}{c}{EFIE} && \multicolumn{3}{c}{CFIE}\\
\cline{3-5} \cline{7-9}
\V && RWG & RWG LU & MR && RWG & RWG LU & MR \\
\hline
\U Satellite at $2$ GHz \subsecref{subsec:emerald} \\ Ram [GB]&& 13.1	& 58.3& 	24.0&&	16.2&	62.3&	26.5\\
\D BiCGStab iterations (threshold $\EE{-3}$) &&  	NC &	NC & NC &&	NC &	14 & 27\\
\hline
\U Satellite $12$ GHz \subsecref{subsec:emerald12} \\ Ram [GB]&& ---	& ---& 	---&&	285& $>$ 650 (extimated)&	477\\
\D BiCGStab iterations (threshold $5\EE{-3}$) &&  	--- &	--- & --- &&	NC &	--- & 99\\
\hline
\U UAV \subsecref{subsec:ucav} \\
Ram [GB]&& 19.4	& 53.2 & 25.6 && 21.5 &	48.1 & 20.1 \\
\D GMRES iterations ( threshold $\EE{-3}$) &&  NC & 538 & 564 &&	NC &	71 & 86\\
\hline
\end{tabular}

\end{center}
\caption{Cases of \secref{sec:realistic}. `NC' means that no convergence is achieved in 1000 iterations, `---' means the case has not been considered.}
\label{tab:realistic}
\end{table*}

\subsection{ESA satellite test case at 2 GHz}
\label{subsec:emerald}

\begin{figure}
\begin{center}
\includegraphics[width=.75\columnwidth]{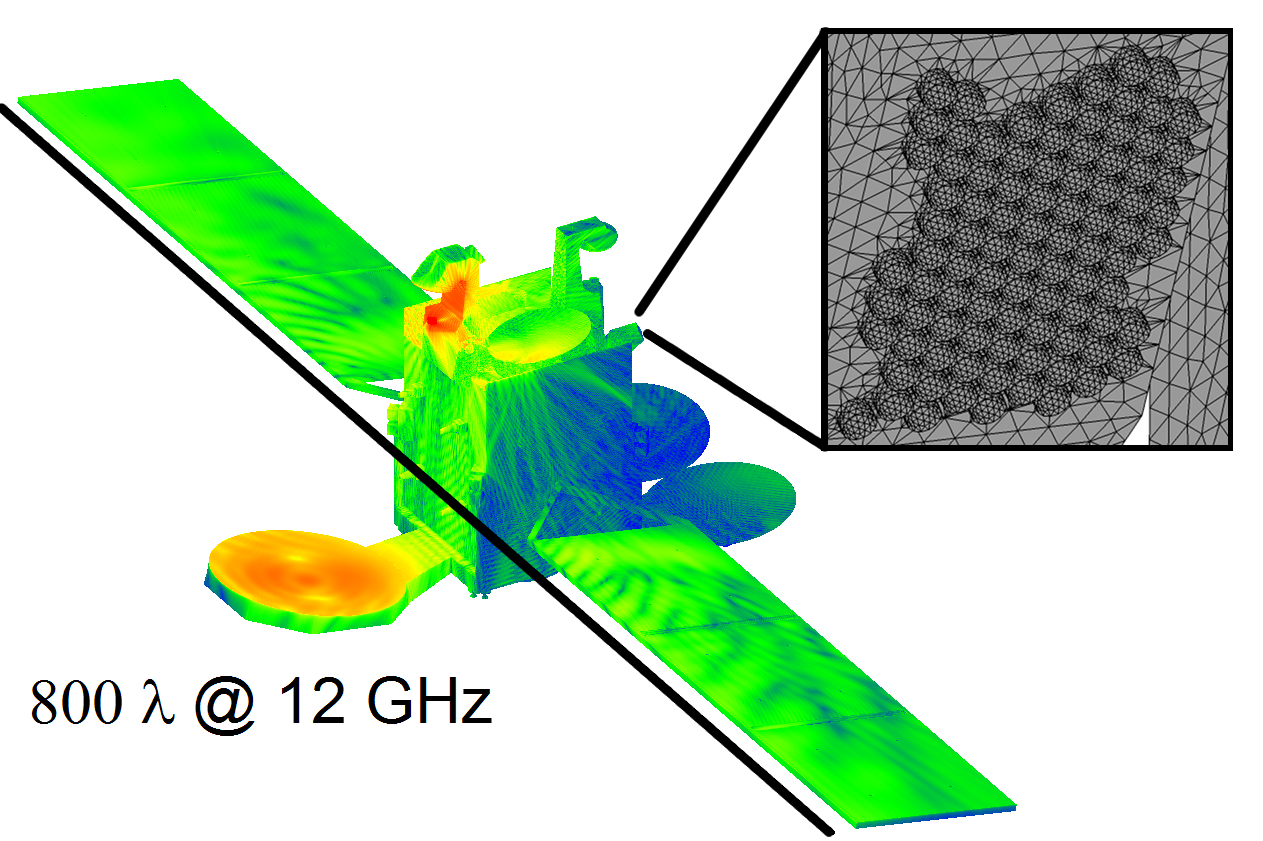}\\
\includegraphics[width=.25\columnwidth]{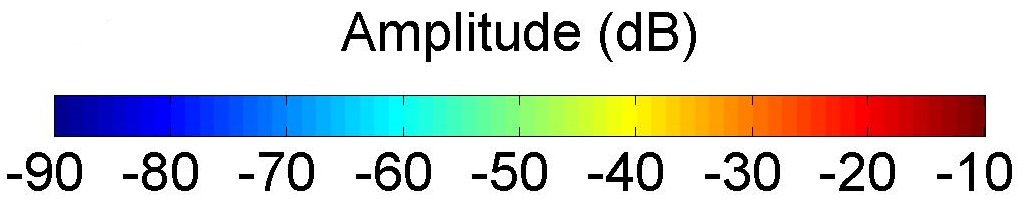}
\end{center}
\caption{Surface current, size, and mesh details for the satellite at $12$ GHz of \subsecref{subsec:emerald12}.}
\label{fig:emerald_current}
\end{figure}

We consider a realistic test case by the European Space Agency, with a helix antenna on top (see \figref{fig:emerald_current}, where mesh and size refer to the $12$ GHz case of \subsecref{subsec:emerald12}). At the frequency of $2$ GHz, the overall dimension is more than $130\lambda$. The mesh has \tocheck{edges from $\lambda/6$ to $\lambda/130$ long}, corresponding to more than $\EE{6}$ unknowns. When enforcing the CFIE, the EFIE is actually enforced on the open parts and on the wire antenna, which amount to nearly $0.2 \%$ of the entire functions. The iterative solver used is the BiCGStab, with a threshold on the residual norm of $\EE{-3}$. \tabref{tab:realistic}, upper part, summarizes the results for different preconditioners and different formulations. The full history of the iterative solver is in \figref{fig:convergence_path}, left panel. In this case, resorting to the CFIE formulation is mandatory to reach convergence, and the MR basis with an algebraic preconditioner on the coarsest functions saves more than half of the memory used in the RWG case, with a negligible increase in the number of iterations needed. 

\subsection{ESA satellite test case at 12 GHz}
\label{subsec:emerald12}

The same structure of \subsecref{subsec:emerald} is considered at $12$ GHz, illuminated with the Spherical Wave Expansion of the feeder of one of the reflectors. The overall size is around $800 \lambda$. The mesh has \tocheck{edges from $\lambda/4$ to $\lambda/114$ long}, giving more than $32\EE{6}$ unknowns. As the EFIE could not reach convergence in the $2$ GHz case of \subsecref{subsec:emerald}, we test the CFIE only. In this very large case, the use MR preconditioning proves the only option to achieve practical convergence, as apparent from \tabref{tab:realistic}; there, we see that we cannot reach convergence in the RWG basis and, on the other hand, the memory required to apply an algebraic preconditioner in the RWG basis is estimated to be more than $650$ GB. Using the MR basis, we can apply the algebraic preconditioner only on the set of coarsest functions and keep the memory requirements feasible with the hardware available, which has $512$ GB of RAM. The obtained surface current is depicted in \figref{fig:emerald_current}.

\subsection{Low profile UAV}
\label{subsec:ucav}
\begin{figure}
\begin{center}
\includegraphics[width=.45\columnwidth]{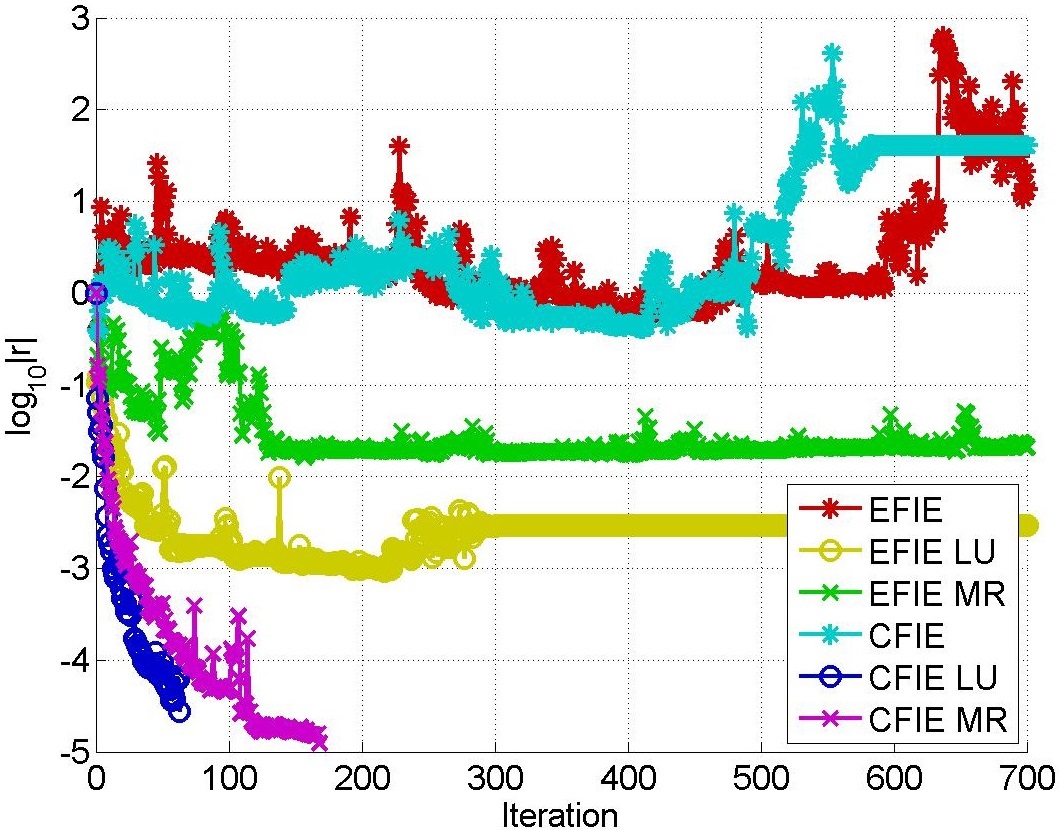}
\includegraphics[width=.45\columnwidth]{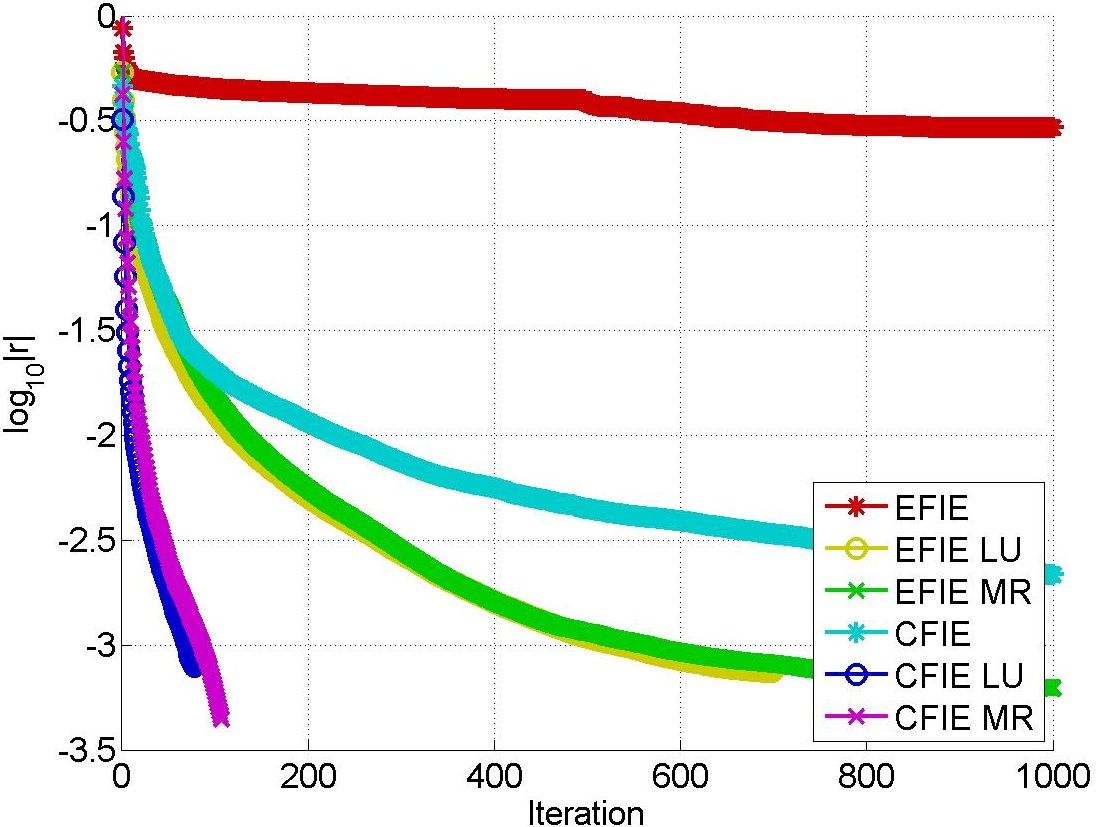}
\end{center}
\caption{Left panel: BiCGStab residual history for the case of \subsecref{subsec:emerald}. Right panel: GMRES residual history for the case of \subsecref{subsec:ucav}. Note that performances should be compared also in terms of RAM memory consumption, see \tabref{tab:realistic}.}
\label{fig:convergence_path}
\end{figure}

\begin{figure}
\begin{center}
\includegraphics[width=.65\columnwidth]{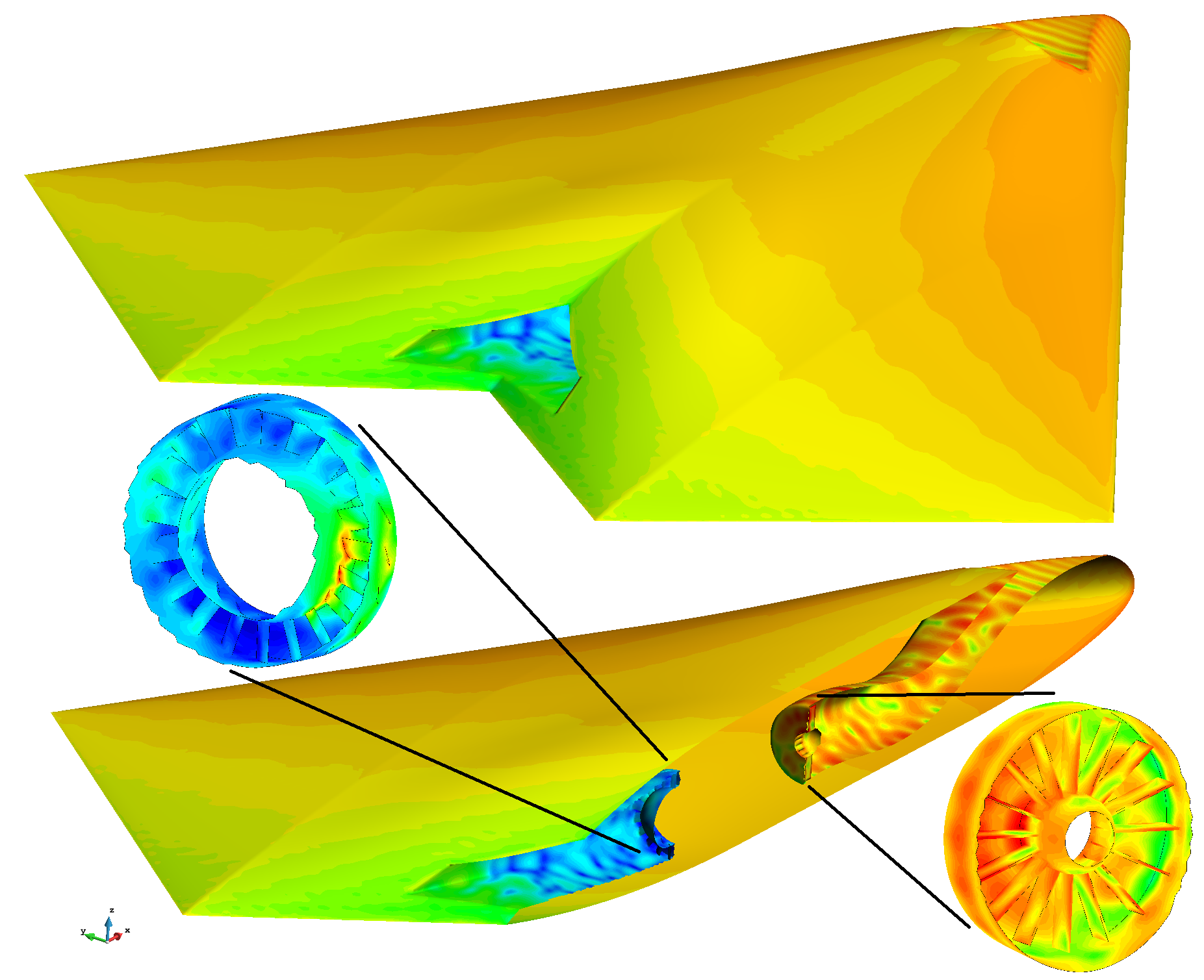}\\
\includegraphics[width=.25\columnwidth]{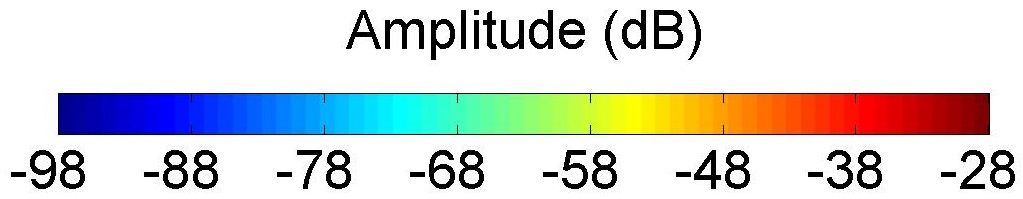}
\end{center}
\caption{Surface current for the UAV of \subsecref{subsec:ucav}, with see-through slice and zooms on the blades of the turbine.}
\label{fig:ucav_current}
\end{figure}

We consider an unmanned aerial vehicle (UAV) with wingspan and tip-to-tail dimension around 30 $\lambda$, illuminated with a linearly polarized plane wave.
The structure is closed, with two deep cavities around the turbine. The considered mesh has edges from $\lambda/7$ to $\lambda/65$ long, resulting in around $7\EE{+5}$ unknowns.  The iterative solver used is the GMRES, with a threshold on the residual norm of $\EE{-3}$. \tabref{tab:realistic}, lower part summarizes the results for different preconditioners and different formulations. \tocheck{The full history of the iterative solver is in \figref{fig:convergence_path}, right panel}. \tocheck{The induced surface current is depicted in \figref{fig:ucav_current}.} In this case, convergence can be achieved with the EFIE formulation as well, with more iterations than with the CFIE. In the MR case, we reach convergence with some more iterations than in the RWG case with a full near-field LU preconditioner, but using less than half the memory. This applies to both the EFIE and the CFIE, meaning that the MR change of basis maintains the good convergence properties of the CFIE and allows an efficient algebraic preconditioner to be applied on a reduced part of the matrix.



\section{Conclusion}
A hierarchical change of basis, developed to improve the convergence of the EFIE for dense discretizations, can be applied in conjunction with an algebraic preconditioner to increase the convergence of the CFIE for electrically large and multiscale structures. Leveraging on the hieararchical change of basis, the preconditioner can be applied on a reduced number of unknowns maintaining the good convergence properties and the immunity to resonance problems of the CFIE. 

We have shown the advantage of the proposed approach on proof-of-principle examples as well as on complex, real-life simulations.

\section*{Acknowledgment}
We thank the European Space Agency for providing the test cases of Sections \ref{subsec:emerald} and \ref{subsec:emerald12}.


\bibliographystyle{IEEEtran}
\bibliography{../../Bibliography}

\end{document}